\documentclass[12pt]{article}
\usepackage{amssymb,amsfonts,amsmath}
\usepackage{graphics}
\usepackage{graphicx}
\usepackage{epsfig}
\usepackage{color}
\usepackage{enumerate}

\newcommand\CC{{\mathbb C}}
\newcommand\GA{\mathcal A}

\newcommand\LL{{\cal L}}

\newcommand\RR{{\mathbb R}}

\newcommand\GX{{\mathcal X}}

\def\beq{\begin{equation}}
\def\eeq{\end{equation}}
\newtheorem{thm}{Theorem}[section]

\newtheorem{lem}[thm]{Lemma}
\newtheorem{cor}[thm]{Corollary}
\newtheorem{rem}[thm]{Remark}

\newcommand\beginpf{\noindent {\bf Proof:} \quad}

\newcommand\re{\mathop{\rm Re}\nolimits}
\newcommand\im{\mathop{\rm Im}\nolimits}

\newcommand\rad{\mathop{\rm Rad}\nolimits}

\def\beginpf{\noindent {\bf Proof :} \quad}
\def\endpf{\rightline{$\square$}}

\newcommand{\ds}{\displaystyle}

\newcommand\inside{\mathop{\rm int}\nolimits}
\renewcommand\phi{\varphi}

\newcommand\Up{\mathcal{U}^{(p)}}
\newcommand\Ep{\mathcal{E}^{(p)}}

\title{Lower estimates near the origin for functional calculus
on operator semigroups}
\author{ I. Chalendar\thanks{I. C. J., UFR de
Math\'ematiques,  Universit\'e  Lyon~1, 43 bld. du 11/11/1918,
69622 Villeurbanne Cedex, France.
  \protect\linebreak[3]
{\tt chalendar@math.univ-lyon1.fr}.},
 J. Esterle\thanks{I.M.B., UMR 5251, Universit\'e de Bordeaux,
351
cours de la Lib\'eration,
33405 Talence Cedex, France.
  \protect\linebreak[3]
{\tt jean.esterle@math.u-bordeaux1.fr}.}
\ and   J.R. Partington\thanks{School of
Mathematics,
University of Leeds, Leeds LS2 9JT, U.K.
\protect\linebreak[3]
{\tt J.R.Partington@leeds.ac.uk}.}
}

\begin{document}

\baselineskip18pt

\maketitle
\begin{abstract}
This paper provides sharp lower estimates near the origin for the functional calculus $F(-uA)$ of a generator $A$ of 
an operator semigroup defined
on the (strictly) positive real line;
here $F$ is given as the Laplace transform of a measure or distribution.
The results are linked to the existence of an identity element or an exhaustive sequence of idempotents
in the Banach algebra generated by the semigroup.
Both the quasinilpotent and non-quasinilpotent
cases are considered, and sharp results are proved extending many in the literature. 
\end{abstract}



\bibliographystyle{plain}

\section{Introduction}

This article is concerned with estimates for $F(-uA)$ where $A$ is the generator of a strongly continuous semigroup $(T(t))_{t>0}$
on a Banach space. Here $F$ is an entire function with $F(0)=0$, given as the Laplace transform of a measure or distribution; the functional calculus defining $F(-uA)$ is given by means of an integral.

This can be seen as providing a wide generalization of results in \cite{BCEP,EM,kalton}, for example, where
quantities such as $\|T(t)-T(2t)\|$ (or its spectral radius) are estimated near the origin. For example, if
$\|T(t)-T(2t)\|<1/4$ on an interval $(0,t_0)$, then, roughly speaking, $(T(t))_{t>0}$ has a bounded infinitesimal
generator (see \cite{BCEP}).

There are two cases to consider, namely, the quasinilpotent and non-quasinilpotent cases, and the techniques used are based
on strong maximum principles for analytic functions.\\

In Section \ref{sec:2}, the case of quasinilpotent semigroups is considered. 
Then
in Section \ref{sec:3} the non-quasinilpotent case is analysed, providing conditions to obtain either an identity
in the closed algebra generated by the semigroup or else an exhaustive sequence $(P_n)_{n \ge 1}$
of idempotents such that $(P_nT(t))_{t}$ has a bounded generator.
Here, the sharpness of the estimates is shown in Remark~\ref{rem:sharp}.\\

\noindent {\bf Notation:}

We write $\CC_+= \{z \in \CC: \re z > 0\}$, and similarly for $\CC_-$.

Let $D(a,R)$ denote the complex disc $\{|z-a| < R\}$.

For a Jordan curve $\Gamma \subset \CC$, we write
$\inside \Gamma$ (the interior of $\Gamma$) for the open set of points in $\CC$ about which the winding number of $\Gamma$ is non-zero.

For $S \subset \CC$ let  $M_c(S)$ denote the space of regular Borel measures having compact support contained in $S$.

\section{Quasinilpotent semigroups}\label{sec:2}

Suppose that $(T(t))_{t>0}$ is a  nontrivial strongly continuous semigroup of quasi\-nilpotent operators acting on a 
Banach space $(\GX,\|.\|)$. Then we write $\GX_0=\left[ \bigcup_{t>0} T(t)\GX\right]^{-\|.\|}$ (closure in norm), and
define a norm
\[
\|x\|_1 = \sup_{t \ge 0} \|T(t) x\|, \qquad \hbox{where} \quad T(0)x=x,
\]
on the subspace 
$ \GX_1:= \{ x \in \GX_0: \|x\|_1 < \infty\}$, which is a Banach space under the norm $\|.\|_1$. 
Further, we write
\beq\label{eq:xtilde}
\widetilde \GX_1:= \left[ \bigcup_{t>0} T(t)\GX \right]^{-\|.\|_1} \subset \GX_1.
\eeq

The following result follows immediately from the main result of \cite{feller53}. It will be used to reduce the
case of a quasinilpotent semigroup to that of a contractive quasinilpotent semigroup.

\begin{thm}\label{thm:feller}
Let $(T(t))_{t>0}$ be a nontrivial strongly continuous semigroup of quasinilpotent operators acting on a 
Banach space $(\GX,\|.\|)$. Then with $(\widetilde \GX_1,\|.\|_1)$ defined as in (\ref{eq:xtilde}) the semigroup
$(T(t)_{| \widetilde \GX_1})_{t>0}$ is a  strongly continuous semigroup of quasinilpotent 
contractions. Moreover for all operators $R$ in the commutant $\{T(t): t>0\}'$ we have
$\|R_{| \widetilde \GX_1}\|_1 \le \|R\|$.
\end{thm}

\subsection{Some complex function theory}

\begin{thm}\label{thm:Jordan}
Let $f: \overline{\CC_+} \to \CC$ be a continuous bounded nonconstant function, holomorphic on $\CC_+$, such that
$f([0,\infty)) \subset \RR$, $f(0)=0$, and with $\lim_{ { x \to \infty},{x \in \RR}} f(x)=0$.

Suppose that $\alpha > 0$ is such that $f(\alpha) \ge |f(x)| $ for all $x \in [0,\infty)$.
Then there exist $a_1,a_2 \in \CC_+$, $a_0 \in (\alpha,a_1)$ and $a_3 \in i\RR$ with $\im a_j > 0$ for  $j=1,2,3$, and $\im a_2=\im a_3$,
and a simple piecewise linear 
Jordan curve $\Gamma_1$ joining $a_1$ to $a_2$ in the upper right half-plane
$\{z \in \CC: \re z>0, \im z >0\}$ and $\delta > 0$ such that

(i) $|f(z)| \ge f(\alpha) + \delta |z-\alpha|^m$ for all $z \in [\alpha,a_1]$, where $m$ (even) is the smallest
positive integer with $f^{(m)}(\alpha) \ne 0$;

(ii) $|f(z)| > |f(a_0)|$ for all $z \in \Gamma_1 \cup [a_2,a_3]$.
\end{thm}

\beginpf
Since $f$ is holomorphic in $\CC_+$, we have, by Taylor's theorem, constants $M>0$ and $\eta>0$ such that
\[
\left | f(z)-f(\alpha)-\frac{(z-\alpha)^m}{m!} f^{(m)}(\alpha) \right | \le M |z-\alpha|^{m+1},
\]
whenever $|z-\alpha| < \eta$. By choosing $a_1$ with $|a_1-\alpha|$ sufficiently small and with argument such that $(a_1-\alpha)^m<0$
(e.g. $\arg(a_1-\alpha)=\pi/m$), we have
condition (i) and hence
$|f(a_1)| > |f(\alpha)|= f(\alpha)$; we may then choose a point $a_0 \in [\alpha,a_1]$ with
\[
|f(\alpha)| < |f(a_0)| < |f(a_1)|.
\]

Let $U= \{z \in \CC_+: |f(z)| > |f(a_0)| \}$. Since $a_1 \in U$, this is a nonempty open set; we let $V$ denote the connected component of $U$
containing $a_1$.

We claim that $\partial V \cap i\RR \ne \emptyset$. Indeed, note that if $z \in \partial V \cap \CC_+$, then $|f(z)| \le |f(a_0)|$, as otherwise if
$|f(z)| > |f(a_0)|$, then by the continuity of $|f|$ a neighbourhood of $z$ is contained in $V$. But we cannot have
$|f(z)| \le |f(a_0)|$ on the whole of $\partial V$, as then by the strong maximum principle 
(see, e.g. \cite[Thm.~9.4]{lieb-loss})  this inequality would hold for all $z \in V$ including $a_1$.

So there exists $a_3 \in \partial V \cap i\RR$ with $\im a_3 > 0$ and $|f(a_3)| > |f(a_0)|$. By the continuity of $|f|$, there exists $\beta>0$ such that
$|f(z)|> |f(a_0)|$ for all $z \in \overline{\CC_+}$ with $|z-a_3| < \beta$. It follows that there is a point $a_2 \in \CC_+$ with
$\im a_2=\im a_3$ and $(a_3,a_2] \subset V$.

Since $V$ is open and connected, it is path-connected, and so we may join $a_1$ to $a_2$ by a polygonal path in $V$.
We may also guarantee that it is simple (does not cross itself): the only difficulty arises if it crosses itself on the arc $(a_0,a_1)$,
when we may replace $a_1$ by the crossing point closest to $a_0$, or if it crosses itself on the line $(a_2,a_3)$, when we
may replace $a_2$ by the crossing point closest to $a_3$.
\endpf

\begin{figure}[htbp]
  \begin{center}
\begin{picture}(0,0)%
\includegraphics{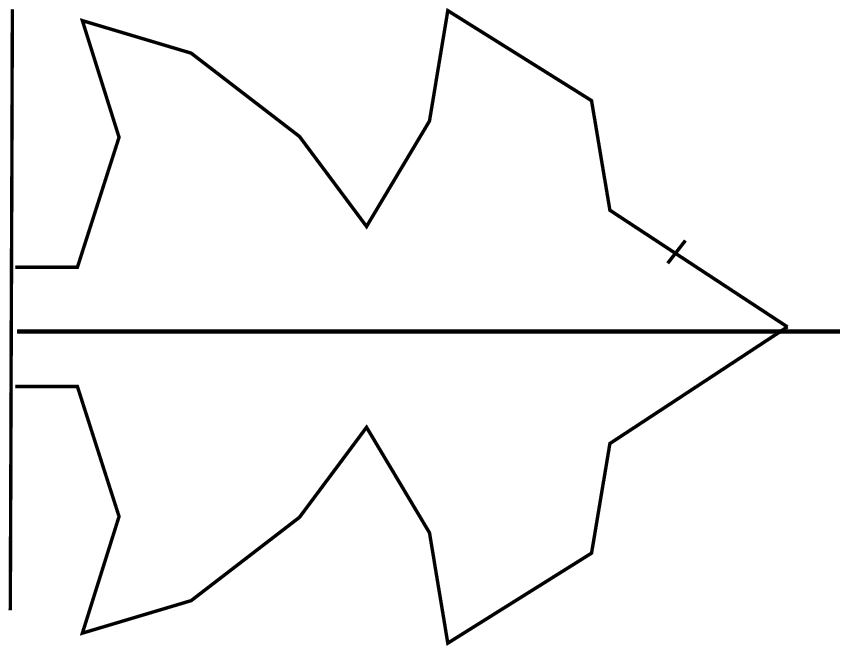}%
\end{picture}%
%
%
\setlength{\unitlength}{7894sp}%
\begingroup\makeatletter\ifx\SetFigFont\undefined%
\gdef\SetFigFont#1#2#3#4#5{%
  \reset@font\fontsize{#1}{#2pt}%
  \fontfamily{#3}\fontseries{#4}\fontshape{#5}%
  \selectfont}%
\fi\endgroup%
\begin{picture}(2157,1545)(83,-723)
\put(2090, 94){\makebox(0,0)[lb]{\smash{{\SetFigFont{14}{16.8}{\rmdefault}{\mddefault}{\updefault}{\color[rgb]{0,0,0}$\alpha$}%
}}}}
\put(1857,276){\makebox(0,0)[lb]{\smash{{\SetFigFont{14}{16.8}{\rmdefault}{\mddefault}{\updefault}{\color[rgb]{0,0,0}$a_0$}%
}}}}
\put(1528,340){\makebox(0,0)[lb]{\smash{{\SetFigFont{14}{16.8}{\rmdefault}{\mddefault}{\updefault}{\color[rgb]{0,0,0}$a_1$}%
}}}}
\put(434,191){\makebox(0,0)[lb]{\smash{{\SetFigFont{14}{16.8}{\rmdefault}{\mddefault}{\updefault}{\color[rgb]{0,0,0}$a_2$}%
}}}}
\put( 98,155){\makebox(0,0)[lb]{\smash{{\SetFigFont{14}{16.8}{\rmdefault}{\mddefault}{\updefault}{\color[rgb]{0,0,0}$a_3$}%
}}}}
\put(1012,449){\makebox(0,0)[lb]{\smash{{\SetFigFont{14}{16.8}{\rmdefault}{\mddefault}{\updefault}{\color[rgb]{0,0,0}$\Gamma_1$}%
}}}}
\end{picture}%
    \caption{The curve constructed in the proof of Theorem \ref{thm:thmhelp}}
    \label{fig:g1}
  \end{center}
\end{figure}

The curve constructed in the proof of Theorem \ref{thm:Jordan} may be seen as the upper part of Figure
\ref{fig:g1}.

We shall also require the following easy result.

\begin{lem}\label{lem:babylem} Let $f \in H(\CC)$ with $f(0)=0$, $f$ nonconstant. Then there is an $r>0$ such that for all $u \in \CC\setminus \{0\}$, with $|u|<r$, we have
\[
\sup_{|z| \le \frac{r}{|u|}} |f(zu)| < \sup_{x \ge 0} |f(xu)|.
\]
\end{lem}
\beginpf
Choose $R>0$ such that $f(z) \ne 0$ for all $z$ with $|z|=R$. Let
$\delta=\inf \{|f(z)|: |z|=R \}$ so that $\delta>0$. Now take $r>0$ such that
$\sup_{|z|\le r} |f(z)| < \delta$, using the continuity of $f$ and the fact that $f(0)=0$.

Then $\sup_{x>0} |f(xu)| \ge \delta$ for all $u \in \CC\setminus \{0\}$, and the conclusion follows.
\endpf

\subsection{The main result}

Recall that if $(T(t))_{t>0}$ is a  uniformly bounded  strongly continuous semigroup with generator $A$, then
\[
(A+\lambda I)^{-1}= - \int_0^\infty e^{\lambda t} T(t) \, dt,
\]
for all $\lambda \in \CC$ with $\re \lambda < 0$.
Here the integral is taken in the sense of Bochner with respect to the strong operator topology.

If, in addition, $(T(t))_{t > 0}$ is quasinilpotent, then
\[
(A+\lambda I)^{-1}= - \int_0^\infty e^{\lambda t} T(t) \, dt,
\]
for all $\lambda \in \CC$.

Similarly, if $\mu \in M_c(0,\infty)$ with Laplace transform
\beq\label{eq:Lmu}
F(s):=\LL \mu(s) = \int_0^\infty e^{-s\xi} \, d\mu(\xi),
\eeq
and $(T(t))_{t>0}$ is a strongly continuous semigroup of bounded operators on $\GX$, then we have a functional
calculus for its generator $A$, defined by
\[
F(-A)= \int_0^\infty T(\xi) \, d\mu(\xi),
\]
in the sense of the strong operator topology; i.e.,
\[
F(-A)x= \int_0^\infty T(\xi)x \, d\mu(\xi), \qquad (x \in \GX),
\]
which exists as a Bochner integral.

\begin{lem}\label{lem:FAFL}
Let $\mu \in M_c(0,\infty)$ and $(T(t))_{t>0}$ a strongly continuous quasinilpotent  semigroup of contractions. Set $F=\LL \mu$. Then we
have for $\re \lambda \ge 0$,
\[
\left\|( F(-A)-F(\lambda)I)(A+\lambda I)^{-1} \right \| \le \int_0^\infty t \, d|\mu| (t).
\]
\end{lem}

\beginpf
We have
\begin{eqnarray*}
F(-A)(A+\lambda I)^{-1} &=& -\left( \int_0^\infty T(t) \, d\mu(t) \right) \left( \int_0^\infty e^{\lambda s}T(s) \, ds \right)\\
&=& -\int_0^\infty  e^{-\lambda t} \left[ \int_0^\infty e^{\lambda(s+t)}T(s+t) \, ds \right] \, d\mu(t) \\
&=& -\int_0^\infty e^{-\lambda t}\left[ \int_t^\infty e^{\lambda v}T(v) \, dv \right] \, d\mu(t),
\end{eqnarray*}
where $v=s+t$. This in turn equals
\begin{eqnarray*}
&& -\int_0^\infty e^{-\lambda t}\left[ \int_0^\infty e^{\lambda v}T(v) \, dv \right] \, d\mu(t)
+\int_0^\infty e^{-\lambda t}\left[ \int_0^t e^{\lambda v}T(v) \, dv \right] \, d\mu(t) \\
&=& F(\lambda)(A+\lambda I)^{-1} + \int_0^\infty \left[ \int_0^t e^{\lambda(v-t)}T(v) \, dv\right] \, d\mu(t).
\end{eqnarray*}
For $\re \lambda \ge 0$, we have
\[
\left\| \int_0^t e^{\lambda(v-t)} T(v) \, dv \right \| \le t,
\]
and so the conclusion follows.

\endpf

The following theorem applies to several examples studied recently in \cite{BCEP,Esterle,EM,kalton}; these include $\mu=\delta_1-\delta_2$, the difference of two Dirac
measures, where $F(s):=\LL \mu(s)=e^{-s}-e^{-2s}$ and $F(-sA)=T(s)-T(2s)$. 
More importantly, the theorem applies to many other examples, such as $d\mu(t)=(\chi_{[1,2]}-\chi_{[2,3]})(t) dt$ and $\mu=\delta_1-3\delta_2+\delta_3+\delta_4$,
which are not accessible with the methods of \cite{BCEP,Esterle,EM,kalton}.

\begin{thm}\label{thm:thmhelp}
Let $\mu \in M_c(0,\infty)$ be a nontrivial real measure such that $\ds \int_0^\infty d\mu(t) = 0$, and let
$(T(t))_{t > 0}$ be a nontrivial strongly continuous quasinilpotent semigroup of bounded operators on a Banach space $\GX$. Set $F=\LL \mu$. Then
there is an $\eta > 0$ such that
\[
\| F(-sA ) \|  > \max_{x \ge 0} |F(x)| \qquad \hbox{for} \quad 0 < s \le \eta.
\]
\end{thm}
 \beginpf
It follows from Theorem \ref{thm:feller} that we may assume without loss of generality that
$(T(t))_{t>0}$ is a strongly continuous  quasinilpotent semigroup of contractions. Let $\alpha>0$ be such that
$|F(x)| \le |F(\alpha)|$ for all $x \ge 0$, and let $s>0$. By considering $-\mu$ instead of $\mu$, if necessary, we may
suppose that $F(\alpha) > 0$.

By Lemma \ref{lem:FAFL} applied to the semigroup $(T(st))_{t>0}$, for $\re\lambda \ge 0$ we obtain 
\[
\left \| F(-sA)(sA+\lambda I)^{-1} \right \| \ge \left\|F(\lambda)(sA+\lambda I)^{-1} \right\| - \int_0^\infty t\, d|\mu|(t).
\]
It follows that 
\[
\|F(-sA)\| \ge |F(\lambda)| - \frac{1}{\left\|(sA+\lambda I)^{-1}\right\|}\int_0^\infty t \, d|\mu|(t)
\]
for $s>0$ and $\re \lambda \ge 0$.

Suppose that there exists $s \in (0,1)$ such that $\|F(-sA)\| \le F(\alpha)$, and consider the simple Jordan curve
\[
\Gamma:=[\alpha,a_1] \cup \Gamma_1 \cup [a_2,a_3] \cup [a_3,\overline{a_3}] \cup [\overline{a_3},\overline{a_2}]
\cup \overline{\Gamma_1} \cup [\overline a_1, \alpha],
\]
where $\Gamma_1,a_1,a_2,a_3$ are defined as in Theorem~\ref{thm:Jordan}, taking $f=F$ (see Figure
\ref{fig:g1}).

We now make various estimates of $\left\|(sA+\lambda I)^{-1}\right\|$ for $\lambda$ on three different parts of $\Gamma$.

1) For $\lambda \in [\alpha,a_1] \cup [\overline{a_1},\alpha]$ we have
\begin{eqnarray*}
F(\alpha) \ge \|F(-sA)\| &\ge& |F(\lambda)| - \frac{1}{\left\|(sA+\lambda I)^{-1}\right\|}\int_0^\infty t \, d|\mu|(t) \\
&\ge& F(\alpha) + \delta|\lambda-\alpha|^m - \frac{1}{\left\|(sA+\lambda I)^{-1}\right\|}\int_0^\infty t \, d|\mu|(t).
\end{eqnarray*}
Hence we obtain
\beq\label{eq:est1}
\left\|(sA+\lambda I)^{-1}\right\| \le \frac{1}{\delta|\lambda-\alpha|^m}\int_0^\infty t \, d|\mu|(t).
\eeq

2) For $\lambda \in \Gamma_1 \cup [a_2,a_3] \cup [\overline{a_3},\overline{a_2}] \cup \overline{\Gamma_1}$ we have
\begin{eqnarray*}
F(\alpha) \ge \|F(-sA)\| &\ge& |F(\lambda)| - \frac{1}{\left\|(sA+\lambda I)^{-1}\right\|}\int_0^\infty t \, d|\mu|(t)\\
&\ge& |F(a_0)| - \frac{1}{\left\|(sA+\lambda I)^{-1}\right\|}\int_0^\infty t \, d|\mu|(t).
\end{eqnarray*}
It follows that
\beq\label{eq:est2}
\left\|(sA+\lambda I)^{-1}\right\| \le \frac{1}{|F(a_0)|-F(\alpha)}\int_0^\infty t \, d|\mu|(t).
\eeq

3) For $x \in \RR$,
\begin{eqnarray*}
\left\| (A+ixI)^{-1} \right\| &=& \left\| -\int_0^\infty T(t) e^{ixt} \, dt \right \| \\
&\le& \int_0^\infty \|T(t)\| \, dt < \infty,
\end{eqnarray*}
since $(T(t))_{t>0}$ is quasinilpotent and contractive.
Therefore 
\beq\label{eq:est3}
\left\|(sA+\lambda I)^{-1}\right\|=\frac{1}{s} \left\| \left(A+\frac{\lambda}{s}I \right)^{-1}\right\| \le \frac{1}{s} \int_0^\infty \|T(t)\| \, dt 
\eeq
for all $\lambda \in [a_3,\overline{a_3}]$.

We can now provide estimates for the quantity $\left\| (\lambda-\alpha)^m \left( A+\frac{\lambda}{s}I \right)^{-1} \right\|$ 
for $\lambda$ on $\Gamma$. Let $R=\max_{\lambda \in \Gamma}|\lambda-\alpha|$.

By (\ref{eq:est1})
\[
\left\| (\lambda-\alpha)^m \left( A+\frac{\lambda}{s}I \right)^{-1} \right\| \le \frac{s}{\delta} \int_0^\infty t \, d|\mu|(t)
\]
for all $\lambda \in [\alpha,a_1] \cup [\overline{a_1},\alpha]$.

By (\ref{eq:est2})
\[
\left\| (\lambda-\alpha)^m \left( A+\frac{\lambda}{s}I \right)^{-1} \right\| \le \frac{sR^m}{|F(a_0)|-F(\alpha)} \int_0^\infty t \, d|\mu|(t)
\]
for all $\lambda \in \Gamma_1 \cup [a_2,a_3] \cup [\overline{a_3},\overline{a_2}] \cup \overline{\Gamma_1}$.

By (\ref{eq:est3})
\[
\left\| (\lambda-\alpha)^m \left( A+\frac{\lambda}{s}I \right)^{-1} \right\| \le R^m\int_0^\infty \|T(t)\| \, dt
\]
for all $\lambda \in [a_3,\overline{a_3}]$.

Since $0 < s \le 1$, for all $z \in \Gamma \cup \inside \Gamma$ we have
\[
\left\|  \left( A+\frac{z}{s}I \right)^{-1} \right\| \le \frac{M}{|z-\alpha|^m},
\]
by the maximum modulus principle,
where
\[
M= \max \left(R^m \int_0^\infty \|T(t)\| \, dt, 
\frac{R^m}{|F(a_0)|-F(\alpha)}\int_0^\infty t \, d|\mu|(t),
\frac{1}{\delta} \int_0^\infty t \, d|\mu|(t)
\right).
\]

Since by hypothesis $F(0)=0$, there is an $r \in (0,\alpha)$ such that
\[
\sup_{|z| \le r} |F(z)| < F(\alpha).
\]
Since $\overline{D(0,r)}\cap \Gamma \cap \CC_+ = \emptyset$, we have $\overline{D(0,r)} \cap \overline{\CC_+}
\subset \Gamma \cup \inside \Gamma$.

Now if $z \in \overline{D(0,r)}$ with $\re z >0$, we have
$|z-\alpha| \ge \alpha-r$, and thus we have
\[
\left\| \left( A + \frac{z}{s}I \right)^{-1} \right\| \le \frac{M}{|z-\alpha|^m} \le \frac{M}{(\alpha-r)^m}.
\]
Also 
\[
\sup_{\re z \le 0} \left\| \left( A + zI \right)^{-1} \right\| \le \int_0^\infty \|T(t)\| \, dt < \infty.
\]
Now, since by Liouville's theorem the function $z \mapsto \left\| \left( A + zI \right)^{-1} \right\| $ is unbounded on $\CC$, it follows
that for all $u>0$ sufficiently small the inequality
\[
\left\|\left (A+\frac{z}{u}I \right)^{-1} \right\| \le \frac{M}{(\alpha-r)^m}
\]
fails to hold for some $z \in \overline{D(0,r)} \cap \CC_+$, depending on $u$.

It follows that there is an $\eta>0$ such that
\[
\|F(-uA)\|> F(\alpha) \qquad \hbox{for all} \quad u \in (0,\eta].
\]
\endpf

If $\mu \in M_c(0,\infty)$ is now a complex measure, then we write $\widetilde F(z)=\overline{F(\overline z)}$, 
which is also an entire function, indeed, the Laplace transform of $\overline{\mu}$.

\begin{cor}\label{cor:26z}
Let $\mu \in M_c(0,\infty)$ be a nontrivial complex measure such that $\ds \int_0^\infty d\mu(t) = 0$, and let
$(T(t))_{t > 0}$ be a nontrivial strongly continuous quasinilpotent semigroup of bounded operators on a Banach space $\GX$. Set $F=\LL \mu$. Then
there is an $\eta > 0$ such that
\[
\| F(-sA )\widetilde F(-sA) \|  > \max_{x \ge 0} |F(x)|^2 \qquad \hbox{for} \quad 0 < s \le \eta.
\]
\end{cor}

\beginpf
The result follows on applying
Theorem~\ref{thm:thmhelp} to the
 real measure $\nu:= \mu * \overline \mu$, whose Laplace transform satisfies
\[
\LL\nu (s) = F(s)\widetilde F(s).
\]
\endpf

We now give similar results for smoother semigroups: let $p>0$ be an integer, and write $\Up$ for the class of semigroups $(T(t))_{t>0}$
such that the mapping $t \mapsto T(t)$ is $p$ times continuous differentiable with respect to the norm topology. Let
$\Ep$ denote the class of distributions of order $p$ with compact support in $(0,\infty)$. For $\phi \in \Ep$ its
action on a $C^p$ function $f$ may be specified in terms of measures $\mu_0,\ldots.\mu_p$, namely,
\[
\langle f, \phi\rangle = \sum_{j=0}^p \int_0^\infty f^{(j)}(t) \, d\mu_j(t).
\]
The Laplace transform of $\phi$ is given by 
\[
F(z): = \LL \phi(z)= \sum_{j=0}^p \int_0^\infty (-z)^j e^{-zt} \, d\mu_j(t).
\]
We write 
$G_j = \LL \mu_j$ and $F_j(z)=(-z)^j G_j(z)$ for each $j$. Likewise
\beq\label{eq:gmeq}
F(-A)= \sum_{j=0}^p A^j G_j(-A) = \sum_{j=0}^p \int_0^\infty A^j T(t) \, d\mu_j(t).
\eeq

We begin with the counterpart of Lemma \ref{lem:FAFL}.

\begin{lem}\label{lem:FAFL2}
Let $p \ge 1$ and $\phi \in \Ep$,  and let $(T(t))_{t>0}$ be a  quasinilpotent  
$\Up$ semigroup of contractions. Set $F=\LL \phi$. Then we
have for $\re \lambda \ge 0$,
\begin{eqnarray*}
\left\|( F(-A)-F(\lambda)I)A^{-p}(A+\lambda I)^{-1} \right \| &\le&\\
 \sum_{m=0}^p c_m\|A^{m-p}\| &+& \sum_{m=0}^p d_m
\left( \sum_{k=0}^{m-1} |\lambda|^k\|A^{m-1-k-p}\| \right),
\end{eqnarray*}
where 
\[
c_m=\int_0^\infty t \, d|\mu_m|(t) \qquad \hbox{and} \quad d_m=\int_0^\infty d|\mu_m|(t)
\]
for $m=0,1,\ldots,p$.
\end{lem}

\beginpf
Write $B:=( F(-A)-F(\lambda)I)A^{-p}(A+\lambda I)^{-1} $. Then by (\ref{eq:gmeq})
we have
\[
B=
 \sum_{m=0}^p A^{-p}(A^m G_m(-A)-(-\lambda)^m G_m(\lambda)I)
 (A+\lambda I)^{-1}.
\]
This can be rewritten as
\begin{eqnarray*}
 && \sum_{m=0}^p A^{m-p}(G_m(-A)-G_m(\lambda))(A+\lambda I)^{-1}\\
&&+ \sum_{m=0}^p G_m(\lambda) A^{-p} (A^m-(-\lambda)^mI)(A+\lambda I)^{-1}.
\end{eqnarray*}
Thus
\begin{eqnarray*}
B= && \sum_{m=0}^p A^{m-p}(G_m(-A)-G_m(\lambda))(A+\lambda I)^{-1} \\
&&+ \sum_{m=0}^p G_m(\lambda) \left[ \sum_{k=0}^{m-1} A^{m-1-k}(-\lambda)^k \right] A^{-p}.
\end{eqnarray*}
Now the first terms can be estimated using Lemma~\ref{lem:FAFL}, and for the second we use the obvious
estimate $|G_m(\lambda)| \le d_m$ for $\re \lambda \ge 0$.
\endpf

\begin{thm}\label{thm:thmhelp2}
Let $p>1$ and $\phi \in \Ep$ be a nontrivial real distribution given by measures 
$\mu_0,\ldots,\mu_p$ such that $\int_0^\infty d\mu_0(t) = 0$, and let
$(T(t))_{t > 0}$ be a  nontrivial
quasinilpotent $\Up$ semigroup of bounded operators on a Banach space $\GX$. Set $F=\LL \phi$. Then
there is an $\eta > 0$ such that
\[
\| F(-sA ) \|  > \max_{x \ge 0} |F(x)| \qquad \hbox{for} \quad 0 < s \le \eta.
\]
\end{thm}

\beginpf
The proof is very similar to the proof of Theorem~\ref{thm:thmhelp}, but using Lemma~\ref{lem:FAFL2}, so we indicate
the changes necessary.
It will be convenient to take $0 < s \le 1$ and to write
\[
K=K(s,\lambda)=  \sum_{m=0}^p c_m\|(sA)^{m-p}\| + \sum_{m=0}^p d_m
\left( \sum_{k=0}^{m-1} |\lambda|^k\|(sA)^{m-1-k-p}\| \right),
\]
noting the dependence on $s$ and $\lambda$.
With the notation of the proof of Theorem~\ref{thm:thmhelp} we have three key estimates: 

1) For $\lambda \in [\alpha,a_1] \cup [\overline{a_1},\alpha]$ we have
\begin{eqnarray*}
F(\alpha) \ge \|F(-sA)\| &\ge& |F(\lambda)| - \frac{K}{\|(sA+\lambda I)^{-1} (sA)^{-p}\|} \\
&\ge& F(\alpha) + \delta|\lambda-\alpha|^m - \frac{K}{\|(sA+\lambda I)^{-1} (sA)^{-p}\|}.
\end{eqnarray*}
Hence we obtain
\beq\label{eq:est1b}
\left\|(sA+\lambda I)^{-1}(sA)^{-p}\right\| \le \frac{K}{\delta|\lambda-\alpha|^m}.
\eeq

2) For $\lambda \in \Gamma_1 \cup [a_2,a_3] \cup [\overline{a_3},\overline{a_2}] \cup \overline{\Gamma_1}$ we have
\begin{eqnarray*}
F(\alpha) \ge \|F(-sA)\| &\ge& |F(\lambda)| - \frac{K}{\left\|(sA+\lambda I)^{-1}(sA)^{-p}\right\|}\\
&\ge& |F(a_0)| -  \frac{K}{\left\|(sA+\lambda I)^{-1}(sA)^{-p}\right\|}.\\
\end{eqnarray*}
It follows that
\beq\label{eq:est2b}
\left\|(sA+\lambda I)^{-1}(sA)^{-p}\right\| \le \frac{K}{|F(a_0)|-F(\alpha)}.
\eeq

3) For $x \in \RR$,
\begin{eqnarray*}
\left\| (A+ixI)^{-1} A^{-p}\right\| &=& \left\| -\int_0^\infty T(t) e^{ixt} \, dt A^{-p}\right \| \\
&\le& \|A^{-1}\|^p \int_0^\infty \|T(t)\| \, dt < \infty,
\end{eqnarray*}
since $(T(t))_{t>0}$ is quasinilpotent and contractive.
Therefore 
\beq\label{eq:est3b}
\left\|(sA+\lambda I)^{-1}(sA)^{-p}\right\| \le \frac{1}{s^{p+1}}\|A^{-1}\|^p \int_0^\infty \|T(t)\| \, dt 
\eeq
for all $\lambda \in [a_3,\overline{a_3}]$.\\

We estimate  $\left\| (\lambda-\alpha)^m \left( A+\frac{\lambda}{s}I \right)^{-1}A^{-p} \right\|$ 
for $\lambda$ on $\Gamma$.

 Let $R=\max_{\lambda \in \Gamma}|\lambda-\alpha|$.
By (\ref{eq:est1b})
\[
\left\| (\lambda-\alpha)^m \left( A+\frac{\lambda}{s}I \right)^{-1} A^{-p}\right\| \le \frac{Ks^{p+1}}{\delta}
\]
for all $\lambda \in [\alpha,a_1] \cup [\overline{a_1},\alpha]$.

By (\ref{eq:est2b})
\[
\left\| (\lambda-\alpha)^m \left( A+\frac{\lambda}{s}I \right)^{-1} A^{-p}\right\| \le \frac{Ks^{p+1}R^m}{|F(a_0)|-F(\alpha)}
\]
for all $\lambda \in \Gamma_1 \cup [a_2,a_3] \cup [\overline{a_3},\overline{a_2}] \cup \overline{\Gamma_1}$.

By (\ref{eq:est3b})
\[
\left\| (\lambda-\alpha)^m \left( A+\frac{\lambda}{s}I \right)^{-1}A^{-p}\right\| \le R^m\|A^{-1}\|^p\int_0^\infty \|T(t)\| \, dt
\]
for all $\lambda \in [a_3,\overline{a_3}]$.

Since $0 < s \le 1$, for all $z \in \Gamma \cup \inside \Gamma$ with $|z|\le r$ we have
\beq\label{eq:goodfact}
\left\|  \left( A+\frac{z}{s}I \right)^{-1} A^{-p}\right\| \le \frac{M}{|z-\alpha|^m} \le \frac{M}{(\alpha-r)^m},
\eeq
by the maximum modulus principle,
where
\[
M= \sup_{\substack{0<s \le 1\\ z \in \Gamma \cup \inside \Gamma}}
 \max \left(\frac{K(s,z)s^{p+1}}{\delta},
\frac{K(s,z)s^{p+1}R^m}{|F(a_0)|-F(\alpha)},
R^m \|A^{-1}\|^p \int_0^\infty \|T(t)\| \, dt
\right),
\]
which is finite.

With this new choice of $M$, the proof is now concluded as for the proof of Theorem~\ref{thm:thmhelp}, using the observation that
$\|(A+zI)^{-1} A^{-p}\|$ is unbounded on $\CC$, and obtaining a contradiction from (\ref{eq:goodfact}).

\endpf

\section{The non-quasinilpotent case}\label{sec:3}

Let $(T(t))_{t>0}$ be a semigroup of non-quasinilpotent operators, and let $\GA_T$ denote the closed (commutative)
algebra generated by the semigroup. We write $\widehat \GA_T$ for the maximal ideal space of $\GA_T$.
Recall that this is compact if and only if $ \GA_T/\rad( \GA_T)$ is unital; otherwise it is locally compact, 
and the function $\widehat a : \chi \mapsto \chi(a)$ is continuous on $\widehat \GA_T$ for every $a \in \GA_T$.



Recall that $\GA_T$ is said to have an {\em exhaustive sequence of idempotents\/} $(P_n)_{n \ge 1}$
if $P_n^2=P_nP_{n+1} = P_n$ for all $n$ and for every $\chi \in \widehat \GA_T$ there is a $p$ such that
$\chi(P_n)=1$ for all $n \ge p$. 

The following result is part of the folklore of the subject, and it is partly contained in \cite[Lem. 3.1]{EM} and \cite[Lem. 3.1]{BCEP}. It enables us to regard $A$ itself
as an element of $C(\widehat \GA_T)$ by defining an appropriate value $\chi(A)=-a_\chi$ for each $\chi \in \widehat \GA_T$.

\begin{lem}\label{lem:em31}
For a strongly continuous and eventually norm-continuous semigroup $(T(t))_{t>0}$ and a nontrivial character $\chi \in \widehat \GA_T$ there is a unique $a_\chi \in \CC$
such that $\chi(T(t))=e^{-ta_\chi}$ for all $t>0$. Moreover, the mapping
$\chi \mapsto a_\chi$ is continuous, and $\chi(F(-uA))=F(ua_\chi)$ in the case that $F=\LL \mu$, as in
(\ref{eq:Lmu}).
\end{lem}

\beginpf
The  existence of $a_\chi$ is given in \cite{EM}, and its uniqueness is clear since the values of $e^{-ta_\chi}$ for $t>0$
determine $a_\chi$ uniquely.

For the continuity, note that 
\[
\chi(T(1))=e^{-a_\chi}
\]
and
\[
\chi\left( e^\lambda\int_1^\infty T(t) e^{-\lambda t} \, dt \right) = \frac{1}{a_\chi+\lambda}e^{-a_\chi}
\]
if $\lambda$ is taken sufficiently large that the integral converges.
Thus if we have a net $\chi_\alpha \to \chi$ then $e^{-a_{\chi_\alpha}} \to e^{-a\chi}$ and 
$ \frac{1}{a_{\chi_\alpha}+\lambda}e^{-a_{\chi_\alpha}} \to \frac{1}{a_\chi+\lambda}e^{-a_\chi}$, which easily
implies that $a_{\chi_\alpha} \to a_\chi$. 

The final observation follows from an easy argument using Bochner integrals.

\endpf

The following result will also be required.

\begin{lem}
\label{lem:shilov}
Let $(T(t))_{t>0}$ be a non-quasinilpotent and eventually norm-con\-tinu\-ous semigroup in a Banach algebra, with infinitesimal generator $A$; let
$\GA_T$ be the subalgebra generated by the semigroup and $\Lambda=\{a_\chi: \chi \in \widehat \GA_T\}$, as
in Lemma \ref{lem:em31}. Then the following conditions are equivalent:
\begin{enumerate}[(i)]
\item $\GA_T$ has an exhaustive sequence of idempotents.
\item For each integer $m \ge 1$ the set $\Lambda_m: = \{ \lambda \in \Lambda: \re \lambda \le m\}$ is contained
in a compact relatively open subset of $\Lambda$.
\end{enumerate}
\end{lem}

\beginpf
$(i)\Rightarrow  (ii):$ Let $\theta: \chi \to a_\chi$ be the homeomorphism given by Lemma~\ref{lem:em31}.
Now $K:=\theta^{-1} (\Lambda_m)$ is a closed subset of
the compact space $\widehat \GA_T \cup \{0\}$, 
since $\re a_\chi \le m$ if and only if $|\chi(T(1))| \ge e^{-m}$.
Hence $K$ is compact. 
 If (i) holds, then for each $\chi\in \widehat \GA_T$, then there is an
$n_\chi>0$ be such that $\chi(P_n)=1$ for $n \ge n_\chi$. So
$\left( \{ \chi \in \widehat \GA_T: \chi(P_n)=1 \} \right)_n$ is an open cover of $K$. By compactness, there is an $N$
such that $\chi(P_N)=1$ for all $\chi \in K$, so $K \subset \{ \chi \in \widehat\GA_T: \chi(P_N)=1\}$, and
$\theta(K)$ is compact, open in $\Lambda$, and contains $\Lambda_m$.\\

$(ii) \Rightarrow (i):$ Conversely, if (ii) holds, then for each $m \ge 1$, the set $\Lambda_m$ is contained in
a compact relatively open set $\Omega_m \subset \Lambda$. So
$\theta^{-1}(\Omega_m) \subset \GA_T$ is compact and open; hence by Shilov's idempotent
theorem \cite[Thm. 2.4.33]{dales00}
there is an idempotent $P_m$ in $\GA_T$ such that
$\chi(P_m)=1$ for $\chi \in \theta^{-1}(\Omega_m)$ and $0$ otherwise. Now $(P_m)_{m \ge 1}$ is an
exhaustive sequence of idempotents in $\GA_T$.

\endpf

\begin{thm}\label{thm:dodgy}
Let $(T(t))_{t>0}$ be a nontrivial strongly continuous and eventually norm-continuous non-quasinilpotent semigroup on a Banach space $\GX$, with generator $A$. Let 
$F=\LL \mu$, where $\mu \in M_c(0,\infty)$ is a real measure such that $\int_0^\infty d\mu=0$. If there exists
$(u_k)_k \subset (0,\infty)$ with $u_k \to 0$ such that
\beq\label{eq:Fsharp}
\rho(F(-u_k A)) < \sup_{x >0} |F(x)|,
\eeq
then the algebra $ \GA_T$ possesses an exhaustive sequence of idempotents $(P_n)_{n \ge 1}$
such that each semigroup $(P_n T(t))_{t>0}$ has a bounded generator.

If, further, $\|F(-u_k A)\| < \sup_{x>0} |F(x)|$, then
$\bigcup_{n \ge 1} P_n \GA_T$ is dense in $\GA_T$.
\end{thm}

\beginpf
For $m \ge 1$ let $\Lambda_m= \{ a_\chi: \re a_\chi \le m\}$, where 
$\chi(T(t))=e^{-a_\chi t}$ as in Lemma \ref{lem:em31}.

Let $K:=\theta^{-1} (\Lambda_m)$, which is a compact set, as seen in the proof of Lemma ~\ref{lem:shilov}.
Hence, $\Lambda_m = \theta(K)$ is compact, since
$\theta: \chi \mapsto a_\chi$ is continuous  by Lemma \ref{lem:em31}.

Therefore, there is an $R_m>0$ such that $\Lambda_m \subset D(0,R_m)$. Note that, by the definition of $\Lambda_m$,
we have
\[
\Lambda \cap \CC_- = \Lambda \cap \CC_- \subset D(0,R_m).
\]
By hypothesis there exists a $u_k > 0$ such that 
\[
|u_k| < \frac{r}{R_m},
\]
where $r>0$ is given by Lemma 
\ref{lem:babylem} and
\[
\rho(F(-u_k A)) < \sup_{x > 0} |F(xu_k)|.
\]
It follows that $\Lambda_m \subset D(0,R_m) \subset D(0,r/|u_k|)$ and 
\[
|F(u_k a_\chi)| < \sup_{x>0} |F(u_k x)|
\]
for all $a_\chi \in \Lambda_m$. Let $\alpha_k$ be such that 
$\sup_{x>0} |F(xu_k)|=|F(\alpha_k u_k)|$.

By Theorem \ref{thm:Jordan} there exists a curve $\Gamma_{k,0}$ in 
$\{z: \re z \ge 0, \im z > 0\}$ joining $\alpha_k \in \RR_+$ to $v_k \in i\RR_+$ with $|v_k|>R_m$
on which $|F(u_k z)| \ge |F(\alpha_k u_k)|$.

\begin{figure}[htbp]
  \begin{center}
\begin{picture}(0,0)%
\includegraphics{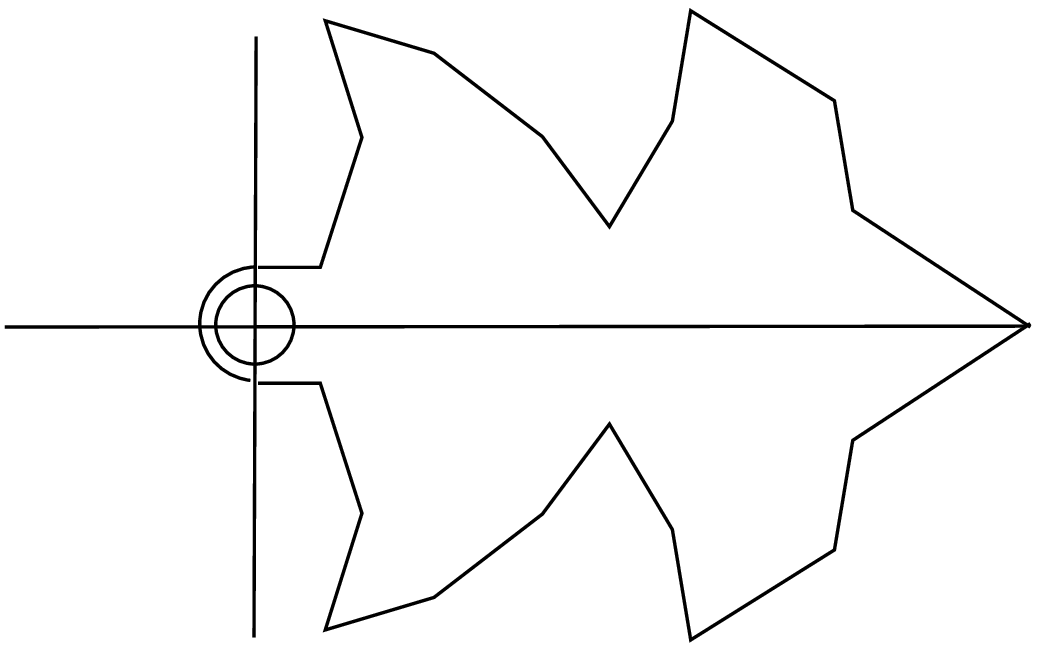}%
\end{picture}%
%
%
\setlength{\unitlength}{7894sp}%
\begingroup\makeatletter\ifx\SetFigFont\undefined%
\gdef\SetFigFont#1#2#3#4#5{%
  \reset@font\fontsize{#1}{#2pt}%
  \fontfamily{#3}\fontseries{#4}\fontshape{#5}%
  \selectfont}%
\fi\endgroup%
\begin{picture}(2485,1534)(439,-730)
\put(2865,111){\makebox(0,0)[lb]{\smash{{\SetFigFont{14}{16.8}{\rmdefault}{\mddefault}{\updefault}{\color[rgb]{0,0,0}$\alpha_k$}%
}}}}
\put(1805,488){\makebox(0,0)[lb]{\smash{{\SetFigFont{14}{16.8}{\rmdefault}{\mddefault}{\updefault}{\color[rgb]{0,0,0}$\Gamma_k$}%
}}}}
\put(934,203){\makebox(0,0)[lb]{\smash{{\SetFigFont{14}{16.8}{\rmdefault}{\mddefault}{\updefault}{\color[rgb]{0,0,0}$v_k$}%
}}}}
\put(1158, 70){\makebox(0,0)[lb]{\smash{{\SetFigFont{14}{16.8}{\rmdefault}{\mddefault}{\updefault}{\color[rgb]{0,0,0}$r/|u_k|$}%
}}}}
\end{picture}%
    \caption{Diagram for the proof of Theorem \ref{thm:dodgy}}
    \label{fig:g2}
  \end{center}
\end{figure}

Let $\Gamma_k= \Gamma_{k,0} \cup \{z \in \CC: \re z<0, \ |z|=|v_k| \} \cup \overline{\Gamma_{k,0}}$ (see
Figure \ref{fig:g2}).

Then $\Lambda \cap \Gamma_k = \emptyset$ since $|F(u_k a_\chi)|=|\chi (F(-u_k A))| < |F(u_k \alpha_k)|$ for $a_\chi \in \Lambda$
and
$|F(u_k z)| > |F(u_k \alpha_k)|$ for $z \in \Gamma_{k,0} \cup \overline{\Gamma_{k,0}}$, so $\Lambda \cap(\Gamma_{k,0} \cup \overline{\Gamma_{k,0}})=\emptyset$.
Also $\Lambda \cap \CC_- = \Lambda_m \cap \CC_-$ so $\Lambda \cap \{z \in \CC_-: |z|=|v_k| \}=\emptyset$.

Now $\Lambda_m = \Lambda \cap \inside \Gamma_k$, which is compact (since $\Lambda \cap \Gamma_k = \emptyset$)
and relatively open in $\Lambda$, so we may now apply Lemma \ref{lem:shilov} to deduce that $\GA_T$ has an exhaustive sequence of idempotents.\\

If $P$ is an idempotent of $\GA_T$, then $\bigcup_{t>0} PT(t) \GA_T$ is dense in the unital
Banach algebra $P \GA_T$. Hence $P \GA_T = \bigcup_{t>0} PT(t) \GA_T$; also $PT(t)$ is invertible in $P \GA_T$ for some, and hence
for all, $t>0$, and then $\lim_{t \to 0+} \|P-PT(t)\|=0$, since the semigroup is eventually continuous.

For the last observation, it follows from Theorem \ref{thm:thmhelp} that $\pi (T(t))=0$ for every $t>0$, where
$\pi: \GA_T \to \GA_T/\bigcup_{n \ge 1}P_n \GA_T$ denotes the canonical surjection.

\endpf

\begin{rem}
A similar result holds for complex measures $\mu \in M_c(0,\infty)$; namely,
we replace
\[
\rho(F(-u_k A)) < \sup_{x >0} |F(x)|,
\]
by the symmetrised version
\[
\rho(F(-u_k A)\widetilde F(-u_k A)) < \sup_{x >0} |F(x)|^2,
\]
as in Corollary~\ref{cor:26z}.
\end{rem}

\begin{rem}\label{rem:sharp}
The following example is given in \cite{BCEP}, and shows that Theorem~\ref{thm:dodgy} is sharp.
Indeed, consider  $C_0[0,1]$, the Banach algebra of all continuous complex-valued functions
on $[0,1]$ that vanish at $0$, equipped with the supremum norm, and 
the semigroup $(S(t))_{t>0}$ defined by $x \mapsto x^t$ for $x \in [0,1]$.

We see that 
\[
F(-uA)= \int_0^\infty S(u\xi)\, d\mu(\xi)= \int_0^\infty x^{u\xi} \, d\mu(\xi),
\]
and for $s>0$ we have
\[
F(s)=\int_0^\infty e^{-s\xi} \, d\mu(\xi).
\]
Thus we obtain equality in \eqref{eq:Fsharp}, choosing $x^u=e^{-s_0}$, where
$s_0$ is defined by
$|F(s_0)|=\sup_{s>0}| F(s)|$.
Note that the norm and spectral radius are equal in $C_0[0,1]$ and that the algebra does not possess any non-trivial
idempotents.
\end{rem}

\section*{Acknowledgements}

This work was partially supported by the ANR project ANR-09-BLAN-0058-01, the
London Mathematical Society (Scheme 2), and the
Institut Camille Jordan. The authors are grateful to the referee for a careful
reading of the paper and some valuable comments.

\end{document}